\documentclass[12pt]{article}
\usepackage{amsmath}
\usepackage{amsfonts}
\usepackage{amssymb}
\usepackage{amscd,epsf, latexsym}
\usepackage[dvips]{graphicx}

\usepackage{amsbsy,eucal, amssymb}
\usepackage{dsfont}
\usepackage{latexsym}

\def\C{{\mathbb{C}}}
\def\R{{\mathbb{R}}}

\newcommand{\andd}{\wedge}
\newcommand{\orr}{\vee}

\newcommand{\iif}{\rightarrow}
\newcommand{\iiff}{\leftrightarrow}
\newtheorem{theorem}{Theorem}
\newtheorem{lemma}[theorem]{Lemma}

\newtheorem{corollary}[theorem]{Corollary}

\newtheorem{examp}{Example}[section]

\newenvironment{definition}{\bigskip\noindent
        \pagebreak[1]\bf Definition\ \rm}{\medskip}

\renewcommand{\phi}{\varphi}

\newcommand{\set}[1]{\{ #1 \}}

\newcommand{\rem}[1]{\relax}

\newcommand{\vvec}{\vec{v}}
\newcommand{\wvec}{\vec{w}}
\newcommand{\zerovec}{\vec{0}}
\newcommand{\abar}{\hat{a}}
\newcommand{\bbar}{\hat{b}}
\newcommand{\cbar}{\hat{c}}
\newcommand{\dbar}{\hat{d}}
\newcommand{\ebar}{\hat{e}}
\newcommand{\fbar}{\hat{f}}
\newcommand{\xbar}{\hat{x}}
\newcommand{\ybar}{\hat{y}}
\newcommand{\zbar}{\hat{z}}

\title
{\Large\bf
\vskip 2truecm
Quantum logic as motivated by quantum computing
}
\author
{\normalsize J. Michael Dunn\footnote{School of Informatics,
Indiana University, Bloomington, IN 47408, dunn@indiana.edu},
\normalsize\; Tobias J. Hagge\footnote{
Department of Mathematics, Indiana University, Bloomington, IN
47405,
thagge@indiana.edu}, \normalsize \;\;
Lawrence S. Mos$\textrm{s}^{\dagger}$\footnote{lsm@cs.indiana.edu}
\normalsize \;\;
and \;Zhenghan Wan$\textrm{g}^{\dagger}$\footnote{zhewang@indiana.edu, partially
supported by NSF EIA 0130388, DMS 0354772, and ARO}
}

\date{}

\begin{document}

\maketitle

\section{Introduction}

Our understanding of Nature comes in layers, so should
the development of logic.  Classic logic is an indispensable part
of our knowledge, and its interactions with computer science have
recently dramatically changed our life.
A new layer of logic has been developing ever
since the discovery of quantum mechanics.  G. D. Birkhoff and
von Neumann introduced quantum logic in a seminal paper in
1936 [BV].
But the definition of quantum logic varies among authors (see [CG]).
How to capture the logic structure inherent in
quantum mechanics is
very interesting and challenging.  Given the close connection between classical logic
and theoretical computer science as exemplified by the coincidence
of computable functions through Turing machines, recursive function theory,
and $\lambda$-calculus, we are interested in how to gain some insights
about quantum logic from quantum computing.  In this note we make
some observations about quantum logic as motivated by quantum computing (see [NC]) and hope
more people will explore this connection.

The quantum logic as envisioned by Birkhoff and von Neumann
is based on the
lattice of closed subspaces of a Hilbert space, usually an
infinite dimensional one. The quantum logic
of a fixed Hilbert space $\mathbb H$ in this note is the variety of all the
true equations with finitely many variables using the connectives
meet, join and negation.
Quantum computing
is theoretically based on quantum systems with finite dimensional Hilbert spaces,
especially the states space of a qubit ${\mathbb C^2}$.  (Actually the qubit
is merely a convenience.  If ${\mathbb C^2}$ is replaced by any
other ${\mathbb C^n}, n\geq 3$, then the same quantum
computing theory will be obtained.)  The $n$-qubit states space
is ${\mathbb
C}^{2^n}$.  It is interesting to understand where the
power of quantum computers could come from.  One possible source is
the exponential growth of the dimensions of the $n$-qubit states space.
Another possibility is the entanglement of quantum states in
${\mathbb C}^{2^n}$ when $n\geq 2$.  Therefore, we ask for a
logical framework to capture the above two sides of quantum computing.
Those questions are only concerned with the static part of quantum
computing.  To include the quantum gates into a logical framework,
we will need a temporal logic. Hence as a first approximation we ask
whether or not the quantum logics of ${\mathbb C}^n$ for all $n$'s are the same,
and what are their relations with
the quantum logic of infinite dimensional Hilbert
spaces.  The quantum logic of $\mathbb C$ reduces to exactly the classical
Boolean logic.  Since the distribution law does not hold in
the quantum logic of ${\mathbb C}^2$, therefore,
 quantum logic of ${\mathbb C}^2$ is different
from that of $\mathbb C$, hence different from classical propositional
logic.  In this note we show that the quantum logic of ${\mathbb
C}^{2^n}$ is always different from that of ${\mathbb
C}^{2^{n+1}}$ for any $n\geq 0$.  We also observe that
quantum logic is not
a finite-valued logic, and
 ${\bf QL}({\mathbb C}^{n})$
is decidable for any $n$.
In the end, we discuss some open problems.

\section{${\bf QL}({\mathbb C}^{2^n}) \neq {\bf QL}({\mathbb C}^{2^{n+1}})$}

Given a Hilbert space $\mathbb H$, let $L_c(\mathbb H)$ be the lattice
of all closed subspaces of $\mathbb H$ with set inclusion as the
partial order relation
$\leq$, and for any two subspaces $p,q \in
L_c(\mathbb H)$, the meet $p\wedge q$ is the set intersection, and
the join $p\vee q$ is the closure of the span of $p\cup q$.  The
closure is necessary in the definition of join when $\mathbb H$ is
infinitely dimensional.  In this case, the span of two subspaces is
not necessarily closed.  For any closed subspace $p$, its negation
$\bar{p}$ is the orthogonal complement.  With the above defined
connectives $\wedge, \vee,\;
\bar{}$ \; on $L_c(\mathbb H)$, $L_c(\mathbb H)$ becomes an orthomodular
lattice.  The maximum and minimum of $L_c(\mathbb H)$ are $\mathbb H$ and $\{0\}$,
respectively.  We will denote them by $\bf 1$ and $\bf 0$.
Recall that an orthomodular lattice is not necessarily
distributive.  It is an ortholattice satisfying the following
orthomodular law:

\begin{equation}
p\wedge [\bar{p}\vee (p\wedge q)]=p\wedge q.
\label{OML}
\end{equation}

The orthomodular law would follow from the distribution law if it holds, so
the orthomodular law is a weakening of the distribution law.
To see the failure of the distribution law $p\vee (q\wedge r)=(p\vee
q)\wedge (p\vee r)$ when $\textrm{dim}(\mathbb H)\geq 2$,
choose $q$ and $r$ to be two distinct lines, and $p$ to be a different
line in the plane spanned by $q,r$.  Then  $(p\vee
q)\wedge (p\vee r)$ is the
plane spanned by $q,r$, while  $p\vee (q\wedge r)$ is just $p$.  Note
also that when $\textrm{dim}(\mathbb H)=1$, then $L_c(\mathbb H)$ is
just $\{\bf{0,1}\}$.

A term $T(p,q,\cdots,r)$ is any formula with finitely many variables
$p,q,\cdots, r$ using connectives $\wedge, \vee,\;  \bar{}$.\;  An
equation $T_1(p,q, \cdots ,r)=T_2(p,q,\cdots, r)$ in $L_c(\mathbb H)$ is
true if the values of the two terms  $T_1$ and $T_2$ always agree
when the variables $p,q,\cdots, r$ range over all $L_c(\mathbb
H)$.

\begin{definition}
Given a Hilbert space $\mathbb H$, the quantum logic ${\bf QL}(\mathbb H)$ is
the variety of all true equations in $L_c(\mathbb H)$.
\end{definition}

In the remaining part of this section, $\mathbb H$ will be always
finitely dimensional.  So it will be ${\mathbb C}^n$ for some $n$,
hence every subspace is closed.

\begin{theorem}
1) $  {\bf QL}({\mathbb C}^2)\neq {\bf QL}({\mathbb C}^4).$

2) For any $i\geq 0$, there exists an equation $\beta \in
{\bf QL}({\mathbb C}^{2^k})$ for any $k\leq i$, and $\beta$ is not in
any ${\bf QL}({\mathbb C}^{2^l})$ if $l>i$.
\label{mainthm}
\end{theorem}

We start with some lemmas.

\begin{lemma}
1) The negation \;$\bar{}$\; is classical, i.e.

a) $ \bar{\bar{p}}=p.$

b) $\overline{p\wedge q}={\bar{p} \vee \bar{q}}.$

c) $\overline{p\vee q}={\bar{p} \wedge \bar{q}}.$

d) $p\wedge {\bar{p}} ={\bf 0}.$

2)  Two subspaces $p=q$ if and only if $ (p\vee q)\wedge
(\bar{p}\vee \bar{q})={\bf 0}$, or equivalently if and only if
$ (\bar{p}\wedge \bar{q})\vee (p\wedge q)={\bf 1}$.
\label{mL1}
\end{lemma}

{\bf Proof:}  The proof of 1) is left to the readers.

To prove 2), denote by $d(s)$ the dimension of a subspace $s$.
If $p=q$, it is a straightforward check.
For the other direction, suppose $(p\vee q)\wedge
(\bar{p}\vee \bar{q})={\bf 0}$.  Note that
 ${\bf 1}=(p\vee q)\vee (\overline{p\vee q})\subseteq
 (p\vee q)\vee (\overline{p\wedge q})=(p\vee q) \vee (\bar{p}\vee
 \bar{q}).$  It follows that ${\bf 1}=(p\vee q) \vee (\bar{p}\vee
 \bar{q}).$
Computing the dimension of $(p\vee q) \vee (\bar{p}\vee \bar{q})$
using the formula $d(s\vee
t)=d(s)+d(t)-d(s\wedge t),$  we obtain
$$n=d(p\vee q)+d(\bar{p}\vee \bar{q})-0=d(p)+d(q)-d(p\wedge
q)+n-d(p\wedge q).$$ Hence  $d(p)+d(q)=2d(p\wedge q)$.  It follows
that $d(p\wedge q)=d(p)=d(q)$.  Since $p\wedge q$ is a subspace of
$p$, and $q$, we have $p\wedge q=p=q$.

\vspace{.2in}

To prove Theorem \ref{mainthm}, we will employ the failure of the distribution
law using Lemma \ref{mL1} (2).  Given three
variables $p,q,r$, we define the distributivity test formula
$\alpha(p,q,r)$ as follows: let
$a=p\vee(q\wedge r)$ and $b=(p\vee q)\wedge (p\vee r)$, and then
define
$$\alpha(p,q,r)=(a\vee b)\wedge (\bar{a}\vee \bar{b}).$$
Note that $a\leq b$, it follows that
$\alpha(p,q,r)=b\wedge \bar{a}=[(p\vee q)\wedge (p\vee r)]\wedge
[\bar{p}\wedge(\bar{q}\vee \bar{r})]\subseteq \bar{p}.$
The distribution law holds if and only if
$\alpha$ is always ${\bf 0}$.  Therefore,
if $\alpha$ does not vanish for some choice of $p,q,r$ in
a Hilbert space, then the distribution law fails for the quantum logic of
this Hilbert space.

\begin{lemma}

Given any three subspaces $p,q,r$ of ${\mathbb C}^n$,
we have ${\textrm{dim}}(\alpha(p,q,r)) \leq \frac{n}{2}.$
\label{mL2}
\end{lemma}
{\bf Proof:}  We know that $\alpha(p,q,r)$ is a subset
of $\bar{p}$.  So if $\textrm{dim}(p)> \frac{n}{2}$, the lemma holds.
Now suppose $\textrm{dim}(p)\leq \frac{n}{2}$.  Since $a\subseteq b$,
so $a$ is perpendicular to $\bar{b}$.  Hence $\bar{\alpha}=a\vee \bar{b}
=a\oplus \bar{b}$.  Writing out
$\bar{\alpha}$, we have $\bar{\alpha}=[\bar{p}\wedge \bar{q})\vee
(\bar{p}\wedge \bar{r})]\oplus [p\vee (q\wedge r)].$   Let
$d(s)$ denote the dimension of any subspace $s$.
A direct computation of the dimension of $\overline{\alpha(p,q,r)}$ results
in the
following formulas:

 $$d(\bar{\alpha})=d(p)+d(q\wedge r)-d(p\wedge q\wedge r)+d(\bar{p}\wedge
\bar{q})+d(\bar{p}\wedge\bar{r})-d(\bar{p}\wedge \bar{q} \wedge
\bar{r})$$
 $$=d(p)+d(q\wedge r)-d(p\wedge q\wedge r)+d(p\vee q\vee
r)+n-d(p\vee q)-d(p\vee r)$$
$$=d(q\wedge r)-d(p\wedge q\wedge r)+d(p\vee q\vee
r)+n-d(p)-d(q)+d(p\wedge q)-d(r)+d(p\wedge r)$$
$$=n-d(p)+d(p\vee q\vee r)-d(q)-d(r)+d(q\wedge r)+d(p\wedge
q)+d(p\wedge r)-d(p\wedge q\wedge r)$$
$$=[n-d(p)]+[d(p\vee q\vee r)-d(q\vee r)]+[d(p\wedge
q)+d(p\wedge r)-d(p\wedge q\wedge r)].$$
The second and third brackets are non-negative numbers.  Since
$d(p)\leq \frac{n}{2}$, so $n-d(p) \geq \frac{n}{2}$.  Hence
$d(\bar{\alpha}) \geq \frac{n}{2}$, and the lemma follows.

\begin{lemma}

For ${\mathbb C}^2$, $\alpha(p,q,r)$ is not ${\bf 0}$ if and only if $p,q,r$
are three distinct lines.
Furthermore, $\alpha(p,q,r)$ is either ${\bf 0}$ or $\bar{p}$.
\label{mL3}
\end{lemma}

{\bf Proof:} Since $q$ and $r$ are symmetric, so we need to check
only the following cases: $p={\bf 0}, {\bf 1}$, $q={\bf 0}, {\bf
1}$, $p=q$, $q=r$, and $p,q,r$ are three distinct lines. Verifications
are left as exercises.  When $\alpha(p,q,r)$ is not $0$, then it
is a one-dimensional subspace of $\bar{p}$.  Since $\bar{p}$ is
also one-dimensional, so $\alpha(p,q,r)$ equals $\bar{p}$ in this
case.

\begin{lemma}
Given two Hilbert spaces $V\subset W$, then ${\bf QL}(V)\supseteq {\bf QL}(W)$.
\label{mL4}
\end{lemma}
{\bf Proof}:  Suppose an equation $T_1=T_2$ of $k$ variables
does not hold for subspaces $p_1,\cdots, p_k$ of $V$.  Then this equation
will not hold for the subspaces $p_i\oplus \bar{V}, i=1,\cdots,k$ of $W$,
where $\bar{V}$ is the orthogonal complement of $V$ in $W$.
This completes the proof.

\vspace{.3in}

\noindent {\bf Proof of Theorem \ref{mainthm}:}

\vspace{.1in}

(1):  By Lemma \ref{mL4}, it suffices to find a true equation in
${\bf QL}({\mathbb C}^2)$, but not in ${\bf QL}({\mathbb C}^4).$
Let
$\beta(p,q,r,s)=\alpha(\alpha(p,q,r),\overline{\alpha(p,q,r)}\wedge
\bar{p},s)$.  We claim $\beta$ is always ${\bf 0}$ in ${\mathbb C}^2$,
but fails for a certain choice of $p,q,r,s$ in ${\mathbb C}^4$.
First we verify that $\beta={\bf 0}$ in ${\mathbb C}^2$.  By Lemma \ref{mL3},
$\alpha(p,q,r)$ is either ${\bf 0}$ or $\bar{p}$.  So is
$\overline{\alpha(p,q,r)}\wedge
\bar{p}$.  It follows from Lemma \ref{mL3} that $\beta={\bf 0}$ since either
at least one of $\alpha(p,q,r)$, $\overline{\alpha(p,q,r)}\wedge
\bar{p}$ is ${\bf 0}$ or they are both $\bar{p}$.  To show $\beta$ is
not always ${\bf 0}$ in ${\mathbb C}^4$, let $\{e_i, i=1,2,3,4\}$ be
an orthonormal basis of ${\mathbb C}^4$ and
$p=\textrm{span of}\;\{ e_1,e_2\}$, $q=\bar{p}$, $r=\textrm{span of}\;\{ e_1,e_2+e_3\}$,
and $s=\textrm{span of}\;\{ e_1,e_3+e_4\}$.  Direct computation shows
$\beta(p,q,r,s)=\textrm{span of}\;\{ e_4\}$, which falsifies $\beta={\bf 0}$
in ${\mathbb C}^4$.

(2)  This general argument will give a different proof for part
(1).  First we explain the restriction of a formula to a new
variable.  Suppose $T$ is a formula for variables
in $L_c({\mathbb C}^n)$, then the restriction of $T$ to a new term $\alpha$,
denoted by
$T|_{\alpha}$, is the following formula:  first using the De Morgan
law, we assume that all negations \;$\bar{}$\; are applied to
single variables, next we replace each variable $p$ and $\bar{p}$
by $p\wedge \alpha$ and $\bar{p}\wedge \alpha$, respectively.
Inductively, we define
$$\alpha^m(p_m,q_m,r_m)=\alpha|_{\alpha^{m-1}}({p_m},{q_m},{r_m}),$$
and $\alpha^1(p_1,q_1,r_1)=\alpha(p_1,q_1,r_1),
\alpha^{m-1}=\alpha^{m-1}(p_{m-1},q_{m-1},r_{m-1}).$
Lemma \ref{mL2} implies that
$$\textrm{dim}(\alpha^m(p_m,q_m,r_m))
\leq \frac{\textrm{dim}(\alpha^{m-1}(p_{m-1},q_{m-1},r_{m-1}))}{2}
\leq \cdots \leq \frac{n}{2^m}$$ if
$\textrm{dim}(\mathbb H)=n$.  In ${\bf QL}({\mathbb C}^{2^i})$,
$\textrm{dim}(\alpha^{i+1})\leq \frac{2^i}{2^{i+1}}<1$, so
$\alpha^{i+1}={\bf 0}$ which gives a true equation in ${\mathbb
C}^{2^i}$.  By Lemma \ref{mL4}, this equation is also true for any $k\leq
i$.  To show it is not true for ${\mathbb C}^{2^{k+1}}$, we notice
that if $p, q, r$ are different subspaces of dimension
$\frac{n}{2}$ and each pair has trivial intersection in ${\mathbb C}^n$, then
$\textrm{dim}\alpha(p,q,r)=\frac{n}{2}$ if $n$ is even.  By choosing
subspaces in ${\mathbb C}^{2^{k+1}}$ this way, we have
$\textrm{dim}(\alpha^{i+1})=\frac{2^{i+1}}{2^{i+1}}=1$.  Now (2)
follows from Lemma \ref{mL4}.

\section{Decidability of ${\bf QL}(\C^n)$}

The modest observation here is that the first-order theory of
each lattice
 $L(\C^n)$ may be reduced to the first-order theory of the reals.
Hence, by Tarski's Theorem [T], we have the decidability of the first-order
theory of each  $L(\C^n)$.  Moreover, we have a stronger result:
the decidability is uniform in $n$.

\begin{theorem}
There is an effective procedure which,
given a natural number $n$ and a sentence
$\phi$
in the  first-order language of complemented lattices,
 gives a sentence
$\phi^*$ in the first-order theory
of $\R$  such that the following are equivalent:
\begin{enumerate}
\item $L(\C^n) \models \phi$.
\item $\R \models \phi^*$.
\end{enumerate}
\label{theorem-main}
\end{theorem}

\begin{corollary}
The first-order theories of the lattices $L(\C^n)$
are decidable (uniformly).
\end{corollary}

We sketch a proof of this result.  The general result would be
messy to write out in full, and so we content ourselves with
a significantly complicated example and some general remarks.

Suppose the we want to know whether or not
\begin{equation}
  L(\C^n)\models (\forall x, y, z) \
\overline{(x \andd y)} \orr z =  y\andd (\bar{z} \orr x).
\label{eq-1}
\end{equation}
The first thing to do is to add new variables $a, b, c, d, e, f$  for
the subterms on both sides, and then write the following:
\begin{equation}
\begin{array}{lcll}  L(\C^n)\models & (\forall x, y, z) & (\forall a,b,c,d,e,f)  \\
& &  [(a = x\andd y)  \andd (b= \bar{a})
\andd (c = b \orr z) \\
& & \andd
(d = \bar{z}) \andd (e = d \orr x) \andd (f = y \andd e)]
 & \iif\ (c = f) \\
\end{array}
\label{eq-2}
\end{equation}
The point of this is that the atomic formulas in (\ref{eq-2}) are
very simple.  In essence, we have taken the complex terms of
(\ref{eq-1}) and ``flattened'' them out.  This move has nothing
to do with the subject at hand, and it is available (and useful)
in other contexts.
In any case, the reader should observe that (\ref{eq-1}) and (\ref{eq-2})
are equivalent.

At this point, we recall that the subspaces of $\C^n$ are the kernels
of $n\times n$ complex matrices.  Our main move in this note is
to replace a quantifier $(\forall x)$
 over subspaces of $\C^n$ with quantification over $n^2$
variables $(\forall x_{11}, \ldots, x_{nn})$
 over $\C$.

In what follows, we shall use the notation $\hat{x}$ for the $(n^2)$-tuple
of variables $x_{11}, \ldots, x_{nn}$.
Further, we shall use the notation $\vvec$ to denote an $n$-tuple
$v_1,\ldots, v_n$ of variables.
We write $\hat{a}\vvec =\zerovec$, for example, to mean the following:
$$  (a_{11} v_1 + \cdots + a_{1n} v_n = 0)
 \andd \cdots \andd
(a_{n1} v_1 + \cdots +a_{nn} v_n = 0)
$$

Now we can render (\ref{eq-2}) as an assertion about $\C$ alone.
\begin{equation}
\begin{array}{lll} \C\models & (\forall \xbar, \ybar, \zbar, \abar,
\bbar,\cbar,\dbar,\ebar,\fbar)  \\
&   [(\forall \vvec)(\abar\vvec =\zerovec \iiff  (\xbar\vvec = \zerovec
\andd \ybar\vvec = \zerovec) ) \\
  & \andd (\forall \vvec)(\bbar\vvec= \zerovec \iiff (\forall \wvec) (\abar\wvec = \zerovec
\iif \vvec\cdot\wvec = 0))\\
 & \andd (\forall \vvec)(\cbar\vvec= \zerovec \iiff
(\exists \wvec^1, \ldots, \wvec^n, r_1, \ldots, r_n) \\
& \qquad
\bigwedge_i(\bbar\wvec^i = \zerovec \orr \zbar\wvec^i =\zerovec)
\andd \bigwedge_{j=1}^n (v_j = \sum_i r_i w^i_j))\\
 & \andd (\forall \vvec)(\dbar\vvec= \zerovec \iiff (\forall \wvec) (\zbar\wvec = \zerovec
\iif \vvec\cdot\wvec = 0))
\\
 & \andd (\forall \vvec)(\ebar\vvec= \zerovec \iiff
(\exists \wvec^1, \ldots, \wvec^n, r_1, \ldots, r_n) \\
& \qquad
\bigwedge_i(\dbar\wvec^i = \zerovec \orr \xbar\wvec^i =\zerovec)
\andd \bigwedge_{j=1}^n (v_j = \sum_i r_i w^i_j))\\
 & \andd (\forall \vvec)(\fbar\vvec =\zerovec \iiff  (\ybar\vvec = \zerovec
\andd \ebar\vvec = \zerovec) )]
\\
  & \iif\ (\forall \vvec) (\cbar \vvec = \zerovec \iiff   \fbar\vvec = \zerovec) \\
\end{array}
\label{eq-3}
\end{equation}
As one can see, the clauses of (\ref{eq-2}) have been replaced by the
more complicated clauses of (\ref{eq-3}).  We explain how this works in
the hardest case, the one for the lattice join.  One of the clauses of (\ref{eq-2})
is $c = b\orr z$.  Recall that this has to do with subspaces of $\C^n$, and we
want to change this to quantification over $(n\times n)$-tuples over $\C$.
One should think of $\cbar$ as the matrix $C = (c_{ij})$.
To say that $c = b\orr z$ is the same as saying that every $\vvec$ such that
$C \vvec = \zerovec$ is in the span of
\begin{equation}
\set{\wvec: B \vvec = \zerovec \mbox{ or } Z \vvec = \zerovec}.
\label{eq-4}
\end{equation}
However, the fact that we are working in $\C^n$ implies that
if $S$ is any set of vectors containing the zero vector, then
  the span of  $S$
of vectors is the same thing as the set of linear combinations of
sets of exactly $n$ vectors from
$S$.  Getting back to our previous point, to say that
 $c = b\orr z$ is the same as saying that every $\vvec$ such that
$C \vvec = \zerovec$ is a linear combination of $n$ vectors from the
set in (\ref{eq-4}).  And this is
$$(\forall \vvec)(\cbar\vvec= \zerovec \iiff
(\exists \wvec^1, \ldots, \wvec^n, r_1, \ldots, r_n)
\bigwedge_{i=1}^n(\bbar\wvec^i = \zerovec \orr \zbar\wvec^i =\zerovec)
\andd \bigwedge_{j=1}^n (v_j = \sum_i r_i w^i_j)).$$
This is exactly the clause which we put into (\ref{eq-3}) on behalf of
 $c = b\orr z$ in (\ref{eq-2}).

The other steps in going from (\ref{eq-2}) to (\ref{eq-3})
are for the lattice meet and for orthogonal complements,
and they are easier.

Now (\ref{eq-3}) is a sentence about $\C$.  But first-order sentences about
the arithmetic properties of $\C$ are reducible to first-order sentence
about $\R$.  Taking all of these observations together, this gives a method
of going from the sentence $\phi$ in (\ref{eq-1}) to a sentence $\phi^*$ in
the theory of $\R$;  $\phi$ and $\phi^*$ have the property
stated in Theorem~\ref{theorem-main}.

Incidentally, our sentence in (\ref{eq-1}) is a universal sentence.
We chose this because the sentences of greatest
interest about the lattice
$L(\C^n)$ are those universal sentences.  But our method also would
work if $\phi$ in  (\ref{eq-1}) had existential quantifiers, or negation.
The details are quite similar to what we have already done, and so we
omit them here.

\section{Open problems}

The state space of $n$-qubits is
the $n$-th tensor power ${({\mathbb C}^2)}^{\otimes n}$ of ${\mathbb C}^2$.
  Quantum computing suggests
the relevance of ${\bf QL}({({\mathbb C}^2})^{\otimes n})$ for all $n$.
 We have the following inclusions of quantum logics:
$$ {\bf QL}(\mathbb C)\supset {\bf QL}({\mathbb C}^2)\supset
{\bf QL}({\mathbb C}^4)\supset\cdots \supset
{\bf QL}({\mathbb C}^{2^n})\supset {\bf QL}({\mathbb C}^{2^{n+1}})\supset\cdots
\supset {\bf QL}({\mathbb C}^{\infty}).$$
We know that the quantum logics of $n$-qubit spaces are
pair-wise distinct.  It is also known that the intersection of
all ${\bf QL}({\mathbb C}^{2^n})$ is
not ${\bf QL}({\mathbb C}^{\infty})$ as it contains the following true equation,
which is one way to define the modular law (see [G]):

\begin{equation}
(p\wedge r) \vee (q\wedge r)=((p\wedge r)\vee q)\wedge r.
\label{ML}
\end{equation}

It is known that the modular law holds
in ${\bf QL}({\mathbb H})$ iff $\mathbb H$ is finitely
dimensional (see [R]), but the orthomodular law (\ref{OML}) does
hold in any Hilbert space.  So we have

$$ \bigcap_{i=0}^{\infty} {\bf QL}({\mathbb C}^{2^i})\supset
{\bf QL}({\mathbb C}^{\infty}).$$

\noindent Some open questions:

\begin{enumerate}

\item  Is ${\bf QL}({\mathbb C}^{\infty})$ decidable?

\item  How to characterize the difference between $ \bigcap_{i=0}^{\infty} {\bf QL}({\mathbb
C}^{2^i})$
and ${\bf QL}({\mathbb C}^{\infty})$?

\item  Are ${\bf QL}({\mathbb C}^n)$ and ${\bf QL}({\mathbb C}^m)$ always
different for $n\neq m$?

\end{enumerate}

It is interesting to characterize the quantum logics of finite
dimensional Hilbert spaces.  Modular lattices are a first
approximation. But there are significant missing ingredients in
the modular lattice formulation of ${\bf QL}({\mathbb C}^n)$.  It has
been shown in [H] that the word problem for modular lattices is not
decidable.  On the other hand, ${\bf QL}({\mathbb C}^n)$ is always decidable.
Another interesting point is the following observation.
A finite set of closed subspaces in
${\mathbb C}^n$ is called {\it a universal test set} if
the truth of any equation is
determined by the evaluations of the subspaces in this set.
It turns out there are no finite universal test sets for ${\bf QL}({\mathbb
C}^m), m\geq 2$.  To see this, consider
the distributivity testing formula $\alpha(p,q,r)$. For simplicity, we
will only give the details for $m=2$.  In order
for the distributivity testing formula $\alpha(p,q,r)$ to
fail, $p, q, r$ must be three distinct lines.
In order for $\alpha(\alpha(\alpha(\alpha(p,q,r),p,s),q,s),r,s)$
to fail, $p, q,r,s$ must
be distinct lines.
Continuing in this manner, we can build a complicated formula
$\gamma$, the failure of which means that the $k$ subspaces
$p,q,\cdots$ are distinct lines.  Since $k$ is arbitrary,
no finite set of lines will
falsify every invalid formula.  This argument works for
any ${\mathbb C}^m, m\geq 2$.


\begin{thebibliography} {[BK]}
\bibitem[CG]{CG}Maria Luisa Dalla Chiara, R. Giuntini, {\it
Quantum logics,} Handbook of Philosophical Logic, vol. {\bf 6},
Chapter 2, Kluwer Academic publisher.

\bibitem[BV]{BV} G. Birkhoff, J. von Neumann,
{\it The logic of quantum mechanics,}  Ann. Math. vol. {\bf 37} 1936,
823-843.

\bibitem[G]{G}G.  Gratzer, Lattice theory: First concepts and
distributive lattices.  W.H. Freeman, 1981.

\bibitem[H]{H}G. Hutchinson,
{\it Recursively unsolvable word problems of modular lattices and diagram-chasing.}
 J. Algebra {\bf 26} (1973), 385--399.

\bibitem[NC]{NC} M. A. Nielsen, I. L. Chuang,  Quantum computation and
quantum information.  Cambridge University Press, Cambridge, 2000.
xxvi+676 pp.

\bibitem[R]{R}M. Redei, Quantum logic in algebraic approach.
Fundamental Theories of Phyiscs, vol 91.  Kluwer Academic
Publishers Gorup, 1998.

\bibitem[T]{T}A. Tarski, A Decision Method for Elementary Algebra and Geometry.
   RAND Corporation, Santa Monica, Calif., 1948.




\end{thebibliography}
\end{document}